\documentclass{article}
\usepackage{amsmath}
\usepackage{amssymb}
\usepackage{amsfonts}

\begin{document}

\title{Superprocesses on ultradistributions}
\date{ }
\author{R. Vilela Mendes\thanks{%
rvmendes@fc.ul.pt; rvilela.mendes@gmail.com} \\
%EndAName
Centro de Matem\'{a}tica e Aplica\c{c}\~{o}es Fundamentais, \\
University of Lisbon, C6 Campo Grande 1749-016 Lisboa, Portugal. }
\maketitle

\begin{abstract}
Stochastic solutions provide new rigorous results for nonlinear PDE's and,
through its local non-grid nature, are a natural tool for parallel
computation. There are two different approaches for the construction of
stochastic solutions: McKean's and superprocesses. In favour of
superprocesses is the fact that they handle arbitrary boundary conditions.
However, when restricted to measures, superprocesses can only be used to
generate solutions for a limited class of nonlinear PDE's. A new class of
superprocesses, namely superprocesses on ultradistributions, is proposed to
extend the stochastic solution approach to a wider class of PDE's.
\end{abstract}

\section{Stochastic solutions and measure-valued processes}

A \textit{stochastic solution} of a linear or nonlinear partial differential
equation is a stochastic process which, starting from a point $x$ in the
domain, generates after time $t$ a boundary measure that, sampling the
initial condition at $t=0$, provides the solution at the point $x$ and time $%
t$. A classical example is the McKean \cite{McKean} construction of a
stochastic solution for the KPP equation,%
\begin{equation}
\frac{\partial v}{\partial t}=\frac{1}{2}\frac{\partial ^{2}v}{\partial x^{2}%
}+v^{2}-v;\hspace{2cm}v\left( 0,x\right) =g\left( x\right) .  \label{1.1}
\end{equation}%
Let $G\left( t,x\right) $ be the Green's operator for the heat equation $%
\partial _{t}v(t,x)=\frac{1}{2}\frac{\partial ^{2}}{\partial x^{2}}v(t,x)$%
\begin{equation*}
G\left( t,x\right) =e^{\frac{1}{2}t\frac{\partial ^{2}}{\partial x^{2}}}
\end{equation*}%
and write the KPP equation in integral form%
\begin{equation}
v\left( t,x\right) =e^{-t}G\left( t,x\right) g\left( x\right)
+\int_{0}^{t}e^{-\left( t-s\right) }G\left( t-s,x\right) v^{2}\left(
s,x\right) ds.  \label{1.2}
\end{equation}%
Denoting by $\left( \xi _{t},P_{x}\right) $ a Brownian motion starting from
time zero and coordinate $x$, Eq.(\ref{1.2}) may be rewritten as%
\begin{eqnarray}
v\left( t,x\right)  &=&\mathbb{E}_{x}\left\{ e^{-t}g\left( \xi _{t}\right)
+\int_{0}^{t}e^{-\left( t-s\right) }v^{2}\left( s,\xi _{t-s}\right)
ds\right\}   \notag \\
&=&\mathbb{E}_{x}\left\{ e^{-t}g\left( \xi _{t}\right)
+\int_{0}^{t}e^{-s}v^{2}\left( t-s,\xi _{s}\right) ds\right\} .  \label{1.3}
\end{eqnarray}%
The \textit{stochastic solution process} is a composite process: a Brownian
motion plus a branching process with exponential holding time $T$, $P\left(
T>t\right) =e^{-t}$. At each branching point the particle splits into two,
the new particles going along independent Brownian paths. At time $t>0$, if
there are $n$ particles located at $x_{1}\left( t\right) ,x_{2}\left(
t\right) ,\cdots x_{n}\left( t\right) $, the solution of (\ref{1.1}) is
obtained by%
\begin{equation}
v\left( t,x\right) =\mathbb{E}_{x}\left\{ g\left( x_{1}(t)\right) g\left(
x_{2}(t)\right) \cdots g\left( x_{n}(t)\right) \right\} .  \label{1.3a}
\end{equation}%
An equivalent interpretation, that corresponds to the second equality in (%
\ref{1.3}), is of a process starting from time $t$ at $x$ and propagating
backwards-in-time to time zero. When it reaches $t=0$ the process samples
the initial condition, that is, it generates a measure $\mu $ at the $t=0$
boundary which yields the solution by (\ref{1.3a}).

The construction of solutions for nonlinear equations, through the
stochastic interpretation of the integral equations, has become an active
field in recent years, applied for example to Navier-Stokes \cite{Jan} \cite%
{Waymire} \cite{Bhatta1} \cite{Ossiander} \cite{Orum}, to Vlasov-Poisson 
\cite{Vilela1} \cite{Vilela2} \cite{Vilela4}, to Euler \cite{Vilela3} to
magnetohydrodynamics \cite{Floriani} and to a fractional version of the KPP
equation \cite{Cipriano}. In addition to providing new exact results for
nonlinear PDE's, the stochastic solutions are also a promising tool for
numerical implementation, in particular for parallel computation using, for
example, the recently developed probabilistic domain decomposition method 
\cite{Acebron1} \cite{Acebron2} \cite{Acebron3}.

There are basically two methods to construct stochastic solutions. The first
method, which will be called the McKean method, illustrated above, is
essentially a probabilistic interpretation of the Picard series. The
differential equations are written as integral equations which are
rearranged in a such a way that the coefficients of the successive terms in
the Picard iteration obey a normalization condition. The Picard iteration is
then interpreted as an evolution and branching process, the stochastic
solution being equivalent to importance sampling of the normalized Picard
series. The second method \cite{Dynkin1} \cite{Dynkin2} \cite{Li} constructs
the boundary measure of a measure-valued stochastic process (a superprocess)
and obtains the solution of the differential equation by a rescaling
procedure. For a detailed comparison of the two methods refer to \cite%
{Vilela5}.

Although being able to handle arbitrary boundary conditions, a limitation of
measure-valued superprocesses is that they can only represent a limited
class of nonlinear partial differential equations\footnote{%
For a detailed account of the nature of the limitations of superprocesses on
measures as related to the positivity of the coefficients in the offspring
generating function see \cite{Vilela5}.}. The main purpose of this paper is
to extend superprocesses from measure-valued to ultradistribution-valued
processes, which lead to a much wider class of stochastic solutions for
partial differential equations.

\section{Ultradistribution-valued superprocesses}

\subsection{Tempered ultradistributions}

A superprocess describes the evolution of a population, without a fixed
number of units, that evolves according to the laws of chance. Given a
countable dense subset $Q$ of $\left[ 0,\infty \right) $ and a countable
dense subset $F$ of a separable metric space $E$, the countable set%
\begin{equation}
M_{1}=\left\{ \sum_{i=1}^{n}\alpha _{i}\delta _{x_{i}}:x_{1}\cdots x_{n}\in
F;\alpha _{1}\cdots \alpha _{n}\in Q;n\geq 1\right\}   \label{S1}
\end{equation}%
is dense (in the topology of weak convergence) on the space $M\left(
E\right) $ of finite Borel measures on $E$ (theorem 1.8 in \cite{Li}). This
is at the basis of the interpretation of the limits of evolving particle
systems as measure-valued superprocesses. On the other hand the
representation of an evolving measure as a collection of measures with point
support is also useful for the construction of solutions of nonlinear
partial differential equations as rescaling limits of measure-valued
superprocesses.

However, as far as representations of solutions of nonlinear PDE's,
superprocesses constructed in the space $M\left( E\right) $ of finite
measures have serious limitations. The set of interaction terms that can be
handled is limited (essentially to $u^{\alpha }\left( x\right) $ with $%
\alpha \leq 2$) and derivative interactions cannot be included as well. The
first obvious generalization would be to construct superprocesses on
distributions of point support, because any such distribution is a finite
sum of deltas and their derivatives \cite{Stein}. However, because in a
general branching process the number of branches is not bounded, one really
needs a framework that can handle arbitrary sums of deltas and their
derivatives. This requirement leads naturally to the space of
ultradistributions of compact support. The space of \textit{%
ultradistributions} $\mathcal{Z}^{\prime }$ is the topological dual of $%
\mathcal{Z}$, a space of test functions for which the Fourier transform is
in $\mathcal{D}$, the space of infinitely differentiable functions of
compact support. The fact that the Fourier transform of $\mathcal{Z}$ has
compact support endows ultradistributions with a rich analytical structure
which makes these "generalized functions" more convenient than distributions
in many applications. An important dense subspace of $\mathcal{Z}^{\prime }$
is the space of \textit{tempered ultradistributions} $\mathcal{U}^{\prime }$
which may be characterized as Fourier transforms of distributions of
exponential type\footnote{%
Distributions which locally are $\mu \left( x\right) =D^{k}\left(
e^{a\left\vert x\right\vert }f\right) $, $f$ bounded and continuous.} \cite%
{Silva1} \cite{Silva2} \cite{Hoskins}.

However, the representation of tempered ultradistributions by analytical
functions is the most convenient one for practical calculations. Let $%
\mathcal{S}$ be the Schwartz space of functions of rapid decrease and $%
\mathcal{U}\subset \mathcal{S}$ those functions in $\mathcal{S}$ that may be
extended into the complex plane as entire functions of rapid decrease on
strips. More precisely 
\begin{equation*}
\mathcal{U}=\cap _{p=0}^{\infty }\mathcal{U}_{p}
\end{equation*}%
with $\mathcal{U}_{p}$ a space of entire functions topologized by the norm $%
\left\Vert \varphi \right\Vert _{p}=\sup_{z\in \Lambda _{p}}\left\{ \left(
1+\left\vert z\right\vert ^{p}\right) \left\vert \varphi \left( z\right)
\right\vert \right\} $, $\Lambda _{p}$ being the open strip $\Lambda
_{p}=\left\{ z\in \mathbb{C}:\left\vert \text{Im}\left( z\right) \right\vert
<p\right\} $. Each $\mathcal{U}_{p}$ space may also be characterized by the
Fourier transform $\mathcal{F}$%
\begin{equation*}
\mathcal{U}_{p}=\left\{ \varphi :\mathcal{F}\left\{ \varphi \right\} \in 
\mathcal{K}_{p}\right\} ,
\end{equation*}%
$\mathcal{K}_{p}$ being the completion of $C^{\infty }$ for the norm $%
\left\Vert \varphi \right\Vert =\max_{0\leq q\leq p}\left\{ \sup \left\vert
e^{p\left\vert x\right\vert }\varphi ^{\left( q\right) }\right\vert \right\} 
$.

$\mathcal{U}^{\prime }$, the topological dual of $\mathcal{U}$, is Silva's
space of tempered ultradistributions \cite{Silva1} \cite{Silva2}. Let $E=%
\mathbb{R}$. Define $B_{\eta }$ as the complement in $\mathbb{C}$ of the
strip $\text{Im}\left( z\right) \leq \eta $%
\begin{equation}
B_{\eta }=\left\{ z:\text{Im}\left( z\right) >\eta \right\}   \label{S2}
\end{equation}%
and $H_{\eta }$ the set of functions which are holomorphic and of polynomial
growth in $B_{\eta }$%
\begin{equation}
\varphi \left( z\right) \in H_{\eta }\Longrightarrow \exists M,\alpha
:\left\vert \varphi \left( z\right) \right\vert <M\left\vert z\right\vert
^{\alpha },\forall z\in B_{\eta }.  \label{S3}
\end{equation}%
Let $H_{\omega }$ be the union of all such spaces%
\begin{equation}
H_{\omega }=\underset{\eta \geq 0}{\cup }H_{\eta }  \label{S4}
\end{equation}%
and in $H_{\omega }$ define the equivalence relation $\Xi $ by%
\begin{equation*}
\varphi \overset{\Xi }{\simeq }\psi \text{ if }\varphi -\psi \text{ is a
polynomial.}
\end{equation*}%
Then, the space of tempered ultradistribution is%
\begin{equation}
\mathcal{U}^{\prime }=H_{\omega }/\Xi   \label{S5}
\end{equation}%
and $\left[ \phi \left( z\right) \right] $ will denote the equivalence
class. The vectorial operations as well as derivation and multiplication by
polynomials, defined on $H_{\omega }$, are compatible with the equivalence
relation and $\mathcal{U}^{\prime }$ becomes a vector space with these
operations.

The Schwartz space $\mathcal{S}^{\prime }$ of tempered distributions may be
identified with a subspace of $\mathcal{U}^{\prime }$ by the Stieltges
transform, that is, a linear mapping of $\mathcal{S}^{\prime }$ on a
subspace $\mathcal{U}^{\prime \ast }$ of $\mathcal{U}^{\prime }$. Namely,
given $\nu \left( x\right) \in \mathcal{S}^{\prime }$%
\begin{equation}
\varphi \left( z\right) =\frac{p\left( z\right) }{2\pi i}\int \frac{\nu
\left( x\right) }{p\left( x\right) \left( x-z\right) }dx+P\left( z\right)
\label{S6}
\end{equation}%
$\left[ \varphi \left( z\right) \right] \in \mathcal{U}^{\prime }$. Here $%
p\left( z\right) $ is a polynomial such that $\nu /p\sim O\left(
x^{-1}\right) $ and $P\left( z\right) $ is an arbitrary polynomial.

Operations on tempered ultradistributions $f\in \mathcal{U}^{\prime }$ are
performed using their analytical images $\varphi \left( z\right) $. For
example $f$ is integrable in $\mathbb{R}$ if there is an $y_{0}\in \mathbb{R}
$ and a $\varphi \left( z\right) $ in $\left[ \varphi \left( z\right) \right]
\in \mathcal{U}^{\prime }$ such that $\varphi \left( x+iy_{0}\right)
-\varphi \left( x-iy_{0}\right) $ is integrable in $\mathbb{R}$ in the sense
of distributions. Then%
\begin{equation}
\left\langle \varphi |g\right\rangle =\oint\limits_{\Gamma _{y_{0}}}\varphi
\left( z\right) g\left( z\right) dz  \label{S7}
\end{equation}%
$\varphi \in \mathcal{U}^{\prime },g\in \mathcal{U}$ and the integral runs
around the boundaries of the strip $\text{Im}\left( z\right) \leq y_{0}$.

An ultradistribution vanishes in an open set $A\in \mathbb{R}$ if $\varphi
\left( x+iy\right) -\varphi \left( x-iy\right) \rightarrow 0$ for $x\in A$
when $y\rightarrow 0$ or, equivalently, if there is an analytical extension
of $\varphi $ to the vertical strip $\text{Re}z\in A$. The support of $\nu $
is the complement in $\mathbb{R}$ of the largest open set where $\nu $
vanishes.

All these notions are easily generalized to $\mathbb{R}^{n}$ \cite{Silva2} 
\cite{Hasumi} by considering products of semiplans as in (\ref{S2}) and the
corresponding polynomial bounds. For the equivalence relation $\Xi $ one
uses pseudopolynomials, that is, functions of the form%
\begin{equation*}
\sum_{j,k}\rho \left( z_{1},\cdots ,\overset{\wedge }{z_{j}},\cdots
,z_{n}\right) z_{j}^{k},
\end{equation*}%
$\overset{\wedge }{z_{j}}$ meaning that this variable is absent from the
arguments of $\rho $.

An ultradistribution $\nu $ in $\mathbb{R}^{n}$ has compact support if there
is a disk $D$ such that any $\varphi $ in $\left[ \varphi \left( z\right) %
\right] \in \mathcal{U}^{\prime }$ has an analytic extension to $\left( 
\mathbb{C}/D\right) ^{n}$. Then the integral in (\ref{S7}) is around a
closed contour containing the support of the ultradistribution.

For the purposes of this paper, the most important property of
ultradistributions of compact support is the fact that any such
ultradistribution has a representation as a series of multipoles \cite%
{Silva2} \cite{Oliveira}%
\begin{equation}
\nu \left( x\right) =\sum_{r_{1}=0}^{\infty }\cdots \sum_{r_{n}=0}^{\infty
}p_{r_{1},\cdots ,r_{n}}\delta ^{\left( r_{1},\cdots ,r_{n}\right) }\left(
x-a\right) .  \label{S8}
\end{equation}%
This result follows from the fact that for compact support one may apply to
the Stieltjes image the Cauchy theorem over a closed contour. The space of 
\textit{tempered ultradistributions of compact support} will be denoted $%
\mathcal{U}_{0}^{\prime }$.

\subsection{Superprocesses}

Let the underlying space of the superprocess be $\mathbb{R}^{n}$. Denote by $%
\left( X_{t},P_{0\,,\nu }\right) $ a branching stochastic process with
values in $\mathcal{U}_{0}^{\prime }$ and transition probability $P_{0,\nu }$
starting from time $0$, $x\in $ $\mathbb{R}^{n}$ and $\nu \in \mathcal{U}%
_{0}^{\prime }$. The process is assumed to satisfy the \textit{branching
property,} that is, given $\nu =\nu _{1}+\nu _{2}$%
\begin{equation}
P_{0,\nu }=P_{0,\nu _{1}}\ast P_{0,\nu _{2}}.  \label{2.5}
\end{equation}%
After the branching $\left( X_{t}^{1},P_{0,\nu _{1}}\right) $ and $\left(
X_{t}^{2},P_{0\,,\nu _{2}}\right) $ are independent and $X_{t}^{1}+X_{t}^{2}$
has the same law as $\left( X_{t},P_{0,\nu }\right) $. In terms of the 
\textit{transition operator }$V_{t}$ operating on functions on $\mathcal{U}$
this would be%
\begin{equation}
\left\langle V_{t}f,\nu _{1}+\nu _{2}\right\rangle =\left\langle V_{t}f,\nu
_{1}\right\rangle +\left\langle V_{t}f,\nu _{2}\right\rangle   \label{2.6}
\end{equation}%
with $V_{t}$ defined by$\ e^{-\left\langle V_{t}f,\nu \right\rangle
}=P_{0,\nu }e^{-\left\langle f,X_{t}\right\rangle }$ or%
\begin{equation}
\left\langle V_{t}f,\nu \right\rangle =-\log P_{0,\nu }e^{-\left\langle
f,X_{t}\right\rangle }  \label{2.7}
\end{equation}%
$f\in \mathcal{U},\nu \in \mathcal{U}_{0}^{\prime }$.

Underlying the usual construction of superprocesses, in the form that is
useful for the representation of solutions of PDE's, there is a stochastic
process with paths that start from a particular point in $\mathbb{R}^{n}$,
then propagate and branch, but the paths preserve the same nature after the
branching. In terms of measures it means that one starts from an initial $%
\delta _{x}$ which at the branching point originates other $\delta ^{\prime
}s$ with at most some scaling factors. It is to preserve this pointwise
interpretation that, in this larger setting, one considers
ultradistributions in $\mathcal{U}_{0}^{\prime }$, because, as stated above,
any ultradistribution in $\mathcal{U}_{0}^{\prime }$ may be represented as a
multipole expansion. Therefore to define the process it suffices to specify
how the branching acts on arbitrary delta derivatives. The construction may
now proceed as in the measure-valued case \cite{Dynkin1} \cite{Dynkin2},
only with a more general branching function.

In $M=\left[ 0,\infty \right) \times \mathbb{R}^{n}$ consider a set $%
Q\subset M$ and the associated exit process $\xi =\left( \xi _{t},\Pi
_{0,x}\right) $ with parameter $k$ defining the lifetime. The process stars
from $x\in \mathbb{R}^{n}$ carrying along an ultradistribution in $\mathcal{U%
}_{0}^{\prime }$ indexed by the path coordinate. At each branching point
(ruled by $\Pi _{0,x}$) of the $\xi _{t}-$process there is a transition
ruled by a $P$ probability in $\mathcal{U}_{0}^{\prime }$ leading to one or
more elements in $\mathcal{U}_{0}^{\prime }$. These $\mathcal{U}_{0}^{\prime
}$ elements are then carried along by the new paths of the $\xi _{t}-$%
process. The whole process stops at the boundary $\partial Q$, finally
defining a exit process $\left( X_{Q},P_{0,\nu }\right) $ on $\mathcal{U}%
_{0}^{\prime }$. If the initial $\nu $ is $\delta _{x}$ and $f\in \mathcal{U}
$ a function on $\partial Q$ one writes%
\begin{equation}
u\left( x\right) =\left\langle V_{Q}f,\delta _{x}\right\rangle =-\log
P_{0,x}e^{-\left\langle f,X_{Q}\right\rangle }  \label{2.10}
\end{equation}%
$\left\langle f,X_{Q}\right\rangle $ being computed on the (space-time)
boundary with the exit ultradistribution generated by the process.

The connection with nonlinear PDE's is established by defining the whole
process to be a $\left( \xi ,\psi \right) -$\textit{superprocess} if $%
u\left( x\right) $ satisfies the equation%
\begin{equation}
u+G_{Q}\psi \left( u\right) =K_{Q}f  \label{2.11}
\end{equation}%
where $G_{Q}$ is the Green operator,%
\begin{equation}
G_{Q}f\left( 0,x\right) =\Pi _{0,x}\int_{0}^{\tau }f\left( s,\xi _{s}\right)
ds  \label{2.12}
\end{equation}%
and $K_{Q}$ the Poisson operator%
\begin{equation}
K_{Q}f\left( x\right) =\Pi _{0,x}1_{\tau <\infty }f\left( \xi _{\tau }\right)
\label{2.13}
\end{equation}%
$\psi \left( u\right) $ means $\psi \left( 0,x;u\left( 0,x\right) \right) $
and $\tau $ is the first exit time from $Q$.

Eq.(\ref{2.11}) is recognized as the integral version of a nonlinear partial
differential equation with the Green operator determined by the linear part
of the equation and $\psi \left( u\right) $ by the nonlinear terms. If the
equation does not possess a natural Poisson clock for the branching (like
the $-v$ term in KPP, Eq.\ref{1.1}) we have to introduce an artificial
lifetime for the particles in the process ($e^{-k}$), which in the end must
vanish ($k\rightarrow \infty $) through a rescaling method.

The superprocess is then constructed as follows: Let $\varphi \left(
s,x;z\right) $ be the branching function at time $s$ and point $x$. Then
denoting $P_{0,x}e^{-\left\langle f,X_{Q}\right\rangle }$ as $e^{-w\left(
0,x\right) }$\ one has%
\begin{equation}
P_{0,x}e^{-\left\langle f,X_{Q}\right\rangle }=e^{-w\left( 0,x\right) }=\Pi
_{0,x}\left[ e^{-k\tau }e^{-f\left( \tau ,\xi _{\tau }\right)
}+\int_{0}^{\tau }dske^{-ks}\varphi \left( s,\xi _{s};e^{-w\left( \tau
-s,\xi _{s}\right) }\right) \right]  \label{2.15}
\end{equation}%
where $\tau $ is the first exit time from $Q$ and $f\left( \tau ,\xi _{\tau
}\right) =\left\langle f,X_{Q}\right\rangle $ is computed with the exit
boundary ultradistribution. Existence of $\left\langle f,X_{Q}\right\rangle $
and hence of $e^{-w\left( 0,x\right) }$ is insured if $f\in \mathcal{U}$ and
the branching function is such that the exit $X_{Q}\in \mathcal{U}%
_{0}^{\prime }$.

For measure-valued superprocesses the branching function would be 
\begin{equation}
\varphi \left( s,y;z\right) =c\sum_{0}^{\infty }p_{n}(s,y)z^{n}  \label{2.14}
\end{equation}%
with $\sum_{n}p_{n}=1$ and $c$ the branching intensity, but now it may be a
more general function.

For the interpretation of the superprocesses as generating solutions of
PDE's, an important role is played by a transformation that uses $%
\int_{0}^{\tau }ke^{-ks}ds=1-e^{-k\tau }$ and the Markov property $\Pi
_{0,x}1_{s<\tau }\Pi _{s,\xi _{s}}=\Pi _{0,x}1_{s<\tau }$, namely%
\begin{eqnarray}
u\left( x,t\right) &=&\Pi _{0,x}\left\{ e^{-kt}u\left( \xi _{t},0\right)
+\int_{0}^{t}ke^{-ks}\Phi \left( \xi _{s},t-s\right) ds\right\}  \notag \\
&=&\Pi _{0,x}\left\{ u\left( \xi _{t},0\right) +k\int_{0}^{t}\left( \Phi
\left( \xi _{s},t-s\right) -u\left( \xi _{s},t-s\right) \right) ds\right\} .
\label{2.15a}
\end{eqnarray}%
Proof of this result is sketched in ch.4 of Ref.\cite{Dynkin1}. A detailed
proof, with the notations used in this paper may be found in \cite{Vilela5}.
Because (\ref{2.15a}) only depends on the Markov properties of the $\left(
\xi _{t},\Pi _{0,x}\right) $ process it also holds in the ultradistribution
context.

Eq.(\ref{2.15a}) converts Eq.(\ref{2.15}) for $e^{-w\left( 0,x\right) }$ into%
\begin{equation}
e^{-w\left( 0,x\right) }=\Pi _{0,x}\left[ e^{-f\left( \tau ,\xi _{\tau
}\right) }+k\int_{0}^{\tau }ds\left[ \varphi \left( s,\xi _{s};e^{-w\left(
\tau -s,\xi _{s}\right) }\right) -e^{-w\left( \tau -s,\xi _{s}\right) }%
\right] \right] .  \label{2.16}
\end{equation}%
Eq.(\ref{2.11}) is now obtained by a limiting process. Let in (\ref{2.16})
replace $w\left( 0,x\right) $ by $\beta w_{\beta }\left( 0,x\right) $ and $f$
by $\beta f$. $\beta $ is interpreted as the mass of the particles and when
the $\mathcal{U}_{0}^{\prime }$-valued process $X_{Q}\rightarrow \beta X_{Q}$
then $P_{\mu }\rightarrow P_{\frac{\mu }{\beta }}$.%
\begin{equation}
e^{-\beta w\left( 0,x\right) }=\Pi _{0,x}\left[ e^{-\beta f\left( \tau ,\xi
_{\tau }\right) }+k_{\beta }\int_{0}^{\tau }ds\left[ \varphi _{\beta }\left(
s,\xi _{s};e^{-\beta w\left( \tau -s,\xi _{s}\right) }\right) -e^{-\beta
w\left( \tau -s,\xi _{s}\right) }\right] \right]  \label{2.17}
\end{equation}%
Two rescaling limits will be used in this paper. The first one, called here
as \textit{type I}, is the one used in the past for superprocesses on
measures, namely it defines%
\begin{equation}
u_{\beta }^{(1)}=\left( 1-e^{-\beta w_{\beta }}\right) /\beta \hspace{0.3cm};%
\hspace{0.3cm}f_{\beta }^{(1)}=\left( 1-e^{-\beta f}\right) /\beta
\label{2.18}
\end{equation}%
and%
\begin{equation}
\psi _{\beta }^{(1)}\left( 0,x;u_{\beta }^{(1)}\right) =\frac{k_{\beta }}{%
\beta }\left( \varphi \left( 0,x;1-\beta u_{\beta }^{(1)}\right) -1+\beta
u_{\beta }^{(1)}\right)  \label{2.19}
\end{equation}%
one obtains from (\ref{2.17})%
\begin{equation}
u_{\beta }^{(1)}\left( 0,x\right) +\Pi _{0,x}\int_{0}^{\tau }ds\psi _{\beta
}^{(1)}\left( s,\xi _{s};u_{\beta }^{(1)}\right) =\Pi _{0,x}f_{\beta
}^{(1)}\left( \tau ,\xi _{\tau }\right)  \label{2.20}
\end{equation}%
that is%
\begin{equation}
u_{\beta }^{(1)}+G_{Q}\psi _{\beta }^{(1)}\left( u_{\beta }^{(1)}\right)
=K_{Q}f_{\beta }^{(1)}.  \label{2.21}
\end{equation}%
One sees from \ref{2.18} that when $\beta \rightarrow 0$, $f_{\beta
}^{(1)}\rightarrow f$ and if $\psi _{\beta }$ goes to a well defined limit $%
\psi $ then $u_{\beta }$ tends to a limit $u$ solution of (\ref{2.11})
associated to a superprocess. Also one sees from (\ref{2.18}) that in the $%
\beta \rightarrow 0$ limit%
\begin{equation}
u_{\beta }^{(1)}\rightarrow w_{\beta }=-\log P_{0,x}e^{-\left\langle
f,X_{Q}\right\rangle }  \label{2.21a}
\end{equation}%
as in Eq.(\ref{2.10}). The superprocess corresponds to a cloud of
ultradistribution "particles" for which both the mass and the lifetime tend
to zero.

An equivalent result is obtained with a \textit{rescaling of type II} 
\begin{equation}
u_{\beta }^{(2)}=\frac{1}{2\beta }\left( e^{\beta w_{\beta }}-e^{-\beta
w_{\beta }}\right) \hspace{0.3cm};\hspace{0.3cm}f_{\beta }^{(2)}=\frac{1}{%
2\beta }\left( e^{\beta f}-e^{-\beta f}\right) .  \label{2.22}
\end{equation}%
Notice that, as before, $u_{\beta }^{(2)}\rightarrow w_{\beta }$ and $%
f_{\beta }^{(2)}\rightarrow f$ when $\beta \rightarrow 0$.

\subsection{Existence of the superprocess}

Existence of the superprocess is existence of a unique solution for the
equation (\ref{2.17}) and its rescaling limit (\ref{2.21}). It will depend
on the appropriate choice of the branching function $\varphi \left(
s,y;z\right) $. For measure-valued processes this function is a polynomial
in $z$, which corresponds to a branching particle system where the offspring
of each particle has the same nature as the parent or, in terms of point
measures, to branching of $\delta $ into other deltas with a positive
coefficient. For ultradistributions of compact support it suffices to
specify the probabilities of branching from an arbitrary delta derivative $%
\delta ^{(n)}$ to other delta derivatives with a positive or negative
coefficient. Because of the multipole representation of ultradistributions
of compact support (Eq. \ref{S8}) any ultradistribution branching may be
obtained by a linear combination of elementary branchings of this type.

Suppose that such a ultradistribution branching is specified. Associated to
the ultradistribution superprocess $\Gamma $ with branching function $%
\varphi $ there is an \textit{enveloping measure superprocess }$\widetilde{%
\Gamma }$ with branching function $\widetilde{\varphi }$ that has the same
branching topology but without any derivative change in the original delta
measure at time zero nor on its sign. General existence conditions for
measure-valued superprocesses have been found in the past \cite{Dynkin3} 
\cite{Fitzsimmons1} \cite{Fitzsimmons2}. Namely $\widetilde{\varphi }$
should have the form%
\begin{equation}
\widetilde{\varphi }\left( s,y:z\right) =-b\left( s,y\right) z-c\left(
s,y\right) z^{2}+\int_{0}^{\infty }\left( e^{-\lambda z}+\lambda z-1\right)
n\left( s,y;d\lambda \right) .  \label{2.23}
\end{equation}%
Suppose that the branching $\widetilde{\varphi }$ for the process $%
\widetilde{\Gamma }$ is of the form (\ref{2.23}). This insures almost sure
existence of $e^{-\left\langle g,\widetilde{X}\right\rangle }$, $\widetilde{X%
}$ being the exit measure generated by the $\widetilde{\Gamma }$ process.
Decomposing the $\widetilde{X}$ measure into the (measure) components
associated to the each one of the delta derivatives of each sign in the
corresponding ultradistribution $\Gamma $ process%
\begin{equation}
\left\langle g,\widetilde{X}\right\rangle =\sum_{n=0}\left( \left\langle g,%
\widetilde{X}_{n}^{\left( +\right) }\right\rangle +\left\langle g,\widetilde{%
X}_{n}^{\left( -\right) }\right\rangle \right) .  \label{2.24}
\end{equation}%
On the other hand in the ultradistribution $\Gamma $ process, the same
computation for the correspondent exit ultradistribution $X$ yields%
\begin{eqnarray}
\left\langle f,X\right\rangle  &=&\sum_{n=0}\left( \left( -1\right)
^{n}\left\langle f^{\left( n\right) },X_{n}^{\left( +\right) }\right\rangle
+\left( -1\right) ^{n-1}\left\langle f^{\left( n\right) },X_{n}^{\left(
-\right) }\right\rangle \right)   \notag \\
&\leq &\sum_{n=0}\left\vert \left\langle f^{\left( n\right) },X_{n}^{\left(
+\right) }\right\rangle \right\vert +\left\vert \left\langle f^{\left(
n\right) },X_{n}^{\left( -\right) }\right\rangle \right\vert   \notag \\
&\leq &M\int_{\partial Q}\sum_{n=0}\left\vert f^{\left( n\right)
}\right\vert .  \label{2.25}
\end{eqnarray}

Hence, one has the following sufficient condition for the existence of a
ultradistribution-valued superprocess:

\textbf{Proposition 1:} A $\mathcal{U}_{0}^{\prime }$
ultradistribution-valued exit superprocess $\Gamma $ exists if the branching
function $\widetilde{\varphi }$ of the associated enveloping exit measure
process $\widetilde{\Gamma }$ is as in Eq. (\ref{2.23}) and the function $f$
is such that the integral over the exit boundary of $\Sigma _{n}\left\vert
f^{(n)}\right\vert $ is finite.

This result imposes some restrictions on the boundary conditions of the
associated nonlinear differential equations, which however are not too
serious. It suffices, for example that $f$ in a finite area boundary $%
\partial Q$ be well approximated by an arbitrary polynomial.

\subsection{Examples: Superprocesses on signed measures and
ultradistributions}

As stated before, because of the multipole expansion property of
ultradistributions of compact support, it suffices to specify how the
branching operates on general delta derivatives.

The variable $z$ that appears in $\varphi _{\beta }\left( s,x;z\right) $ is
in fact $z=e^{-\beta w\left( \tau -s,\xi _{s}\right)
}=P_{0,x}e^{-\left\langle \beta f,X\right\rangle }$. When restricting the
superprocess to measures, the delta measure, at each branching point, may at
most branch into other deltas (with positive coefficients) and therefore $%
\varphi \left( s,x;z\right) $ must be a sum of monomials in $z$, with
positive coefficients to have a probability interpretation. When one
generalizes to $\mathcal{U}_{0}^{\prime }$, changes of sign and transitions
from deltas to their derivatives are allowed. In the end, the exponential $%
e^{-\left\langle \beta f,X\right\rangle }$ will be computed by evaluation of
the function on the ultradistribution that reaches the boundary. To find out
the equation that is represented by the process one then computes $\psi
_{\beta }\left( 0,x;u_{\beta }\right) $ of Eq.(\ref{2.19}) for the
corresponding $\varphi \left( s,x;z\right) $ in the $\beta \rightarrow 0$
limit. Recalling that $\varphi \left( s,x;z\right) =\varphi _{\beta }\left(
s,\xi _{s};e^{-\beta w\left( \tau -s,\xi _{s}\right) }\right) $ and $%
z=e^{-\beta w_{\beta }}$, one concludes that there are basically two new
transitions at the branching points:

1) A change of sign in the point support ultradistribution%
\begin{equation}
e^{\left\langle \beta f,\delta _{x}\right\rangle }=e^{\beta f\left( x\right)
}\rightarrow e^{\left\langle \beta f,-\delta _{x}\right\rangle }=e^{-\beta
f\left( x\right) }  \label{3.1}
\end{equation}%
which corresponds to%
\begin{equation}
z\rightarrow \frac{1}{z}  \label{3.1a}
\end{equation}%
and

2) A change from $\delta ^{(n)}$ to $\pm \delta ^{(n+1)}$, for example%
\begin{equation}
e^{\left\langle \beta f,\delta _{x}\right\rangle }=e^{\beta f\left( x\right)
}\rightarrow e^{\left\langle \beta f,\pm \delta _{x}^{\prime }\right\rangle
}=e^{\mp \beta f^{\prime }\left( x\right) }  \label{3.2}
\end{equation}%
which corresponds to%
\begin{equation}
z\rightarrow e^{\mp \partial _{x}\log z}.  \label{3.2a}
\end{equation}%
Case 1) corresponds to an extension of superprocesses on measures to
superprocesses on signed measures and the second to superprocesses in $%
\mathcal{U}_{0}^{\prime }$. Another possible transformation would be one
decreasing the order of the derivatives in the $\delta $'s. This might be
useful to generate solutions of integrodifferential equations, but will not
be dealt with here.

Now, referring back to Eq.(\ref{2.11}), one knows that to obtain a
superprocess that generates solutions of a particular nonlinear partial
differential equation amounts to finding a branching function $\varphi
\left( 0,x;z\right) $ which, in the scaling limit, generates a $\psi \left(
0,x;u\right) $ identical to the nonlinear term of the equation. How this
provides stochastic representations of solutions for a larger class of
PDE's, is illustrated by two results:

\textbf{Proposition 2:} The superprocess with branching function%
\begin{equation}
\varphi \left( 0,x;z\right) =p_{1}e^{\partial _{x}\log z}+p_{2}e^{-\partial
_{x}\log z}+p_{3}z^{2}  \label{3.3}
\end{equation}%
provides a solution to the equation%
\begin{equation}
\frac{\partial u}{\partial t}=\frac{1}{2}\frac{\partial ^{2}u}{\partial x^{2}%
}-2u^{2}-\frac{1}{2}\left( \partial _{x}u\right) ^{2}  \label{3.3b}
\end{equation}%
whenever the boundary function $\left. u\right\vert _{\partial Q}$ satisfies
the condition of proposition 1.

Proof: The branching function $\varphi \left( 0,x;z\right) $ means that at
each branching point, with probability $p_{1}$ a derivative is added to the
propagating ultradistribution, with probability $p_{2}$ a derivative is
added plus a change of sign and with probability $p_{3}$ the
ultradistribution branches into two identical ones. The branching function $%
\widetilde{\varphi }$ of the associated enveloping measure process is $%
\left( p_{1}+p_{2}\right) z+p_{3}z^{2}$, therefore belonging to the class of
branchings in Eq. (\ref{2.23}).

Using now the transformation and rescaling (\ref{2.18}) one has, for small $%
\beta $%
\begin{equation}
z\rightarrow e^{\mp \partial _{x}\log z}=e^{\mp \partial _{x}\log \left(
1-\beta u_{\beta }^{(1)}\right) }=1\pm \beta \partial _{x}u_{\beta }^{(1)}+%
\frac{\beta ^{2}}{2}\left\{ \left( \partial _{x}u_{\beta }^{(1)}\right)
^{2}\pm \partial _{x}u_{\beta }^{(1)2}\right\} +O\left( \beta ^{3}\right)
\label{3.4a}
\end{equation}%
\begin{equation}
z\rightarrow z^{2}=\left( 1-\beta u_{\beta }^{(1)}\right) ^{2}=1-2\beta
u_{\beta }^{(1)}+\beta ^{2}u_{\beta }^{(1)2}.  \label{3.4b}
\end{equation}%
Then, computing $\psi _{\beta }\left( 0,x;u_{\beta }^{(1)}\right) $ with $%
p_{1}=p_{2}=\frac{1}{4}$ and $p_{3}=\frac{1}{2}$ one obtains%
\begin{eqnarray}
\psi _{\beta }\left( 0,x;u_{\beta }^{(1)}\right) &=&\frac{k_{\beta }}{\beta }%
\left( \varphi ^{(1)}\left( 0,x;z\right) -z\right)  \notag \\
&=&\frac{k_{\beta }}{\beta }\left( \varphi ^{(1)}\left( 0,x;1-\beta u_{\beta
}^{(1)}\right) -1+\beta u_{\beta }^{(1)}\right)  \notag \\
&=&\frac{k_{\beta }}{\beta }\left( \frac{1}{8}\beta ^{2}\left( \partial
_{x}u_{\beta }^{(1)}\right) ^{2}+\frac{1}{2}\beta ^{2}u_{\beta
}^{(1)2}+O\left( \beta ^{3}\right) \right)  \label{3.5}
\end{eqnarray}%
meaning that, with $k_{\beta }=\frac{4}{\beta }$, the superprocess provides,
in the $\beta \rightarrow 0$ limit, a solution to the equation (\ref{3.3b})

\textbf{Proposition 3: }The superprocess associated to the branching function%
\begin{equation}
\varphi \left( 0,x;z\right) =p_{1}z^{2}+p_{2}\frac{1}{z}  \label{3.5a}
\end{equation}%
provides a solution to the equation%
\begin{equation}
\frac{\partial u}{\partial t}=\frac{1}{2}\frac{\partial ^{2}u}{\partial x^{2}%
}+u^{3}  \label{3.5b}
\end{equation}%
whenever the boundary function $\left. u\right\vert _{\partial Q}$ satisfies
the condition of proposition 1.

Proof: The branching function $\varphi $ means that with probability $p_{1}$
the ultradistribution branches into two identical ones and with probability $%
p_{2}$ it changes its sign. Therefore, in this case, one is simply extending
the superprocess construction to signed measures. The $\widetilde{\varphi }$
branching function is $p_{2}z+p_{1}z^{2}$.

Here the rescaling of Eq.(\ref{2.22}) is used. With $z=e^{\beta w_{\beta }}$
one has%
\begin{eqnarray}
z &=&-2\beta u_{\beta }^{(2)}+2\sqrt{\beta ^{2}u_{\beta }^{(2)2}+1}  \notag
\\
&=&2-2\beta u_{\beta }^{(2)}+\beta ^{2}u_{\beta }^{(2)2}+O\left( \beta
^{4}\right)   \label{3.8a}
\end{eqnarray}%
and%
\begin{eqnarray}
\frac{1}{z} &=&2\beta u_{\beta }^{(2)}+2\sqrt{\beta ^{2}u_{\beta }^{(2)2}+1}
\notag \\
&=&2+2\beta u_{\beta }^{(2)}+\beta ^{2}u_{\beta }^{(2)2}+O\left( \beta
^{4}\right) .  \label{3.8b}
\end{eqnarray}%
For the integral equation, instead of (\ref{2.20}), one has%
\begin{equation}
u_{\beta }^{(2)}\left( 0,x\right) +\Pi _{0,x}\int_{0}^{\tau }ds\psi _{\beta
}^{(2)}\left( s,\xi _{s};u_{\beta }^{(2)}\right) =\Pi _{0,x}f_{\beta
}^{(2)}\left( \tau ,\xi _{\tau }\right)   \label{3.9}
\end{equation}%
with%
\begin{equation}
\psi _{\beta }^{(2)}\left( 0,x;u_{\beta }^{(2)}\right) =k_{\beta }\left( 
\frac{1}{2\beta }\left( \varphi \left( 0,x;z\right) -\varphi \left( 0,x;%
\frac{1}{z}\right) \right) -u_{\beta }^{(2)}\right) .  \label{3.10}
\end{equation}%
Let now the branching function $\varphi \left( 0,x;z\right) $ be as stated
in (\ref{3.5a}) 
\begin{equation*}
\varphi \left( 0,x;z\right) =p_{1}z^{2}+p_{2}\frac{1}{z}.
\end{equation*}%
Using (\ref{3.8a}) and (\ref{3.8b}) one computes $\psi _{\beta }^{(2)}\left(
0,x;u_{\beta }^{(2)}\right) $ obtaining%
\begin{equation}
\psi _{\beta }^{(2)}\left( 0,x;u_{\beta }^{(2)}\right) =k_{\beta }\left\{
-p_{1}8u_{\beta }^{(2)}\left( 1+\frac{1}{2}\beta ^{2}u_{\beta
}^{(2)2}\right) +p_{2}u_{\beta }^{(2)}-u_{\beta }^{(2)}+O\left( \beta
^{4}\right) \right\}   \label{3.12}
\end{equation}%
and with $p_{1}=\frac{1}{10};p_{2}=\frac{9}{10}$ and $k_{\beta }=\frac{5}{%
2\beta ^{2}}$ one obtains in the in the $\beta \rightarrow 0$ limit%
\begin{equation}
\psi _{\beta }^{(2)}\left( 0,x;u_{\beta }^{(2)}\right) \rightarrow -u_{\beta
}^{(2)3}  \label{3.13}
\end{equation}%
meaning that this superprocess provides a solution to the equation (\ref%
{3.5b})

In conclusion: Extending the superprocess construction to signed measures
and ultradistributions, stochastic solutions are obtained for a much larger
class of partial differential equations.

\end{document}